
\newcount\secno
\newcount\prmno
\newif\ifnotfound
\newif\iffound

\def\section#1{\vskip0,8truecm
               \global\def\currenvir{section}
               \global\advance\secno by1\global\prmno=0
               {\bf \number\secno. {#1}}
              }

\def\subsection{\global\def\currenvir{subsection}
                \global\advance\prmno by1
                \smallskip \ind{(\number\secno.\number\prmno) }}
\def\subsec{\global\def\currenvir{subsection}
                \global\advance\prmno by1
                { (\number\secno.\number\prmno)\ }}
\def\numo{ (\number\secno.\number\prmno.}
\def\proclaim#1{\global\advance\prmno by 1
                {\bf #1 \the\secno.\the\prmno$.-$ }}

\long\def\th#1 \enonce#2\endth{%
   \medbreak\proclaim{#1}{\it #2}\global\def\currenvir{th}\smallskip}

\def\bib#1{\rm #1}
\long\def\thr#1\bib#2\enonce#3\endth{%
\medbreak{\global\advance\prmno by 1\bf#1\the\secno.\the\prmno\ 
\bib{#2}$\!\!.-$ } {\it
#3}\global\def\currenvir{th}\smallskip}
\def\rem#1{\global\advance\prmno by 1
{\it #1} \the\secno.\the\prmno$.-$ }


\magnification 1250
\pretolerance=500 \tolerance=1000  \brokenpenalty=5000
\mathcode`A="7041 \mathcode`B="7042 \mathcode`C="7043
\mathcode`D="7044 \mathcode`E="7045 \mathcode`F="7046
\mathcode`G="7047 \mathcode`H="7048 \mathcode`I="7049
\mathcode`J="704A \mathcode`K="704B \mathcode`L="704C
\mathcode`M="704D \mathcode`N="704E \mathcode`O="704F
\mathcode`P="7050 \mathcode`Q="7051 \mathcode`R="7052
\mathcode`S="7053 \mathcode`T="7054 \mathcode`U="7055
\mathcode`V="7056 \mathcode`W="7057 \mathcode`X="7058
\mathcode`Y="7059 \mathcode`Z="705A
\def\spacedmath#1{\def\packedmath##1${\bgroup\mathsurround =0pt##1
\egroup$}
\mathsurround#1
\everymath={\packedmath}\everydisplay={\mathsurround=0pt}}
\def\nospacedmath{\mathsurround=0pt
\everymath={}\everydisplay={} } \spacedmath{2pt}
\def\qfl#1{\buildrel {#1}\over {\longrightarrow}}
\def\phfl#1#2{\normalbaselines{\baselineskip=0pt
\lineskip=10truept\lineskiplimit=1truept}\nospacedmath\smash 
{\mathop{\hbox to 8truemm{\rightarrowfill}}
\limits^{\scriptstyle#1}_{\scriptstyle#2}}}
\def\hfl#1#2{\normalbaselines{\baselineskip=0truept
\lineskip=10truept\lineskiplimit=1truept}\nospacedmath
\smash{\mathop{\hbox to
12truemm{\rightarrowfill}}\limits^{\scriptstyle#1}_{\scriptstyle#2}}}

\def\mono{\lhook\joinrel\mathrel{\longrightarrow}}
\def\iso{\vbox{\hbox to .8cm{\hfill{$\scriptstyle\sim$}\hfill}
\nointerlineskip\hbox to .8cm{{\hfill$\longrightarrow $\hfill}} }}

\def\sdir_#1^#2{\mathrel{\mathop{\kern0pt\oplus}\limits_{#1}^{#2}}}
\def\pprod_#1^#2{\raise
2pt \hbox{$\mathrel{\scriptstyle\mathop{\kern0pt\prod}\limits_{#1}^{#2}}$}}

\font\eightrm=cmr8         \font\eighti=cmmi8
\font\eightsy=cmsy8        \font\eightbf=cmbx8
\font\eighttt=cmtt8        \font\eightit=cmti8
\font\eightsl=cmsl8        \font\sixrm=cmr6
\font\sixi=cmmi6           \font\sixsy=cmsy6
\font\sixbf=cmbx6\catcode`\@=11
\def\eightpoint{%
  \textfont0=\eightrm \scriptfont0=\sixrm \scriptscriptfont0=\fiverm
  \def\rm{\fam\z@\eightrm}%
  \textfont1=\eighti  \scriptfont1=\sixi  \scriptscriptfont1=\fivei
  \def\oldstyle{\fam\@ne\eighti}\let\old=\oldstyle
  \textfont2=\eightsy \scriptfont2=\sixsy \scriptscriptfont2=\fivesy
  \textfont\itfam=\eightit
  \def\it{\fam\itfam\eightit}%
  \textfont\slfam=\eightsl
  \def\sl{\fam\slfam\eightsl}%
  \textfont\bffam=\eightbf \scriptfont\bffam=\sixbf
  \scriptscriptfont\bffam=\fivebf
  \def\bf{\fam\bffam\eightbf}%
  \textfont\ttfam=\eighttt
  \def\tt{\fam\ttfam\eighttt}%
  \abovedisplayskip=9pt plus 3pt minus 9pt
  \belowdisplayskip=\abovedisplayskip
  \abovedisplayshortskip=0pt plus 3pt
  \belowdisplayshortskip=3pt plus 3pt 
  \smallskipamount=2pt plus 1pt minus 1pt
  \medskipamount=4pt plus 2pt minus 1pt
  \bigskipamount=9pt plus 3pt minus 3pt
  \normalbaselineskip=9pt
  \setbox\strutbox=\hbox{\vrule height7pt depth2pt width0pt}%
  \normalbaselines\rm}\catcode`\@=12

\newcount\noteno
\noteno=0
\def\up#1{\raise 1ex\hbox{\sevenrm#1}}
\def\note#1{\global\advance\noteno by1
\footnote{\parindent0.4cm\up{\number\noteno}\
}{\vtop{\eightpoint\baselineskip12pt\hsize15.5truecm\noindent
#1}}\parindent 0cm}
\font\san=cmssdc10
\def\ext{\hbox{\san \char3}}

\def\pc#1{\tenrm#1\sevenrm}
\def\tx{\kern-1.5pt -}
\def\cqfd{\kern 2truemm\unskip\penalty 500\vrule height 4pt depth 0pt width
4pt\medbreak} 
\def\virg{\raise
.4ex\hbox{,}}
\def\decale#1{\smallbreak\hskip 28pt\llap{#1}\kern 5pt}
\def\no{n\up{o}\kern 2pt}
\def\ind{\par\hskip 0.8truecm\relax}
\def\indp{\par\hskip 0.5truecm\relax}

\def\rond{\kern 1pt{\scriptstyle\circ}\kern 1pt}

\def\det{\mathop{\rm det}\nolimits}
\def\Pic{\mathop{\rm Pic}\nolimits}

\def\dim{\mathop{\rm dim}\nolimits}

\def\Sing{\mathop{\rm Sing}\nolimits}
\frenchspacing
\input amssym.def
\vsize = 25.2truecm
\hsize = 16truecm
\voffset = -.5truecm
\parindent=0cm
\baselineskip15pt
\overfullrule=0pt
\def\P{\Bbb P}
\def\F{\Bbb F}
\def\C{\Bbb C}
\def\Z{\Bbb Z}

\def\mr{{\cal M}_r}
\def\jg{J^{g-1}}

\centerline{\bf  Vector bundles on curves and  theta functions}
\smallskip
\smallskip \centerline{Arnaud {\pc BEAUVILLE}} 
\vskip1.2cm

{\bf Introduction}
\smallskip
\ind These notes survey the relation between the moduli spaces of 
vector bundles on a curve $C$ and the spaces of (classical) theta functions
on the Jacobian $J$ of $C$. The connection appears when one tries to
describe the moduli space
$\mr$ of rank $r$ vector bundles with trivial determinant as a projective
variety in an explicit way (as opposed to the somewhat non-constructive way
provided by  GIT). The Picard group of the moduli space is generated by one
line bundle, the {\it determinant line bundle}
${\cal L}$; thus the natural maps from $\mr$ to projective space are
those defined by the linear systems $|{\cal L}^k|$, and in the first instance the
map
$\varphi_{\cal L}:\mr\dasharrow |{\cal L}|^*$.  The key point is that this map
can be identified with the {\it theta map} 
$$\theta:\mr\dasharrow |r\Theta|$$ 
which associates to a general bundle $E\in\mr$ its {\it theta divisor} $\Theta
_E$, an element of the linear system $|r\Theta |$ on $J$ -- we will recall the
precise definitions below. This description turns out to be sufficiently
manageable to get some  information on the behaviour of this map, at least
when $r$ or $g$ are small. 
 We will describe the results which have been obtained so far -- most
of them fairly recently. Thus these notes 
can be viewed as a sequel to [B2], though with a more precise
focus on the theta map.  For the convenience of the reader we have
made this paper independent of [B2], by   recalling in \S 1 the necessary
definitions. Then we discuss the indeterminacy locus of $\theta $ (\S 2),
the case
$r=2$ (\S 3), the case $g=2$ (\S 4), and the higher rank case (\S 5). Finally, as
in [B2] we will propose a small list of questions and conjectures related to the
topic (\S 6).

\font\cag=cmbsy10
\font\grit=cmmib7 
\section{The moduli space $\hbox{\cag\char77}_{\hbox{\grit\char114}}$ and
the theta map}
\subsection Throughout this paper $C$ will be a complex curve of genus $g\ge
2$. We denote by $J$ its Jacobian variety, and by $J^{k}$ the variety
(isomorphic to $J$) parametrizing line bundles of degree $k$
on $C$. 
\ind For $r\ge 2$, we denote by $\mr$ the moduli space of semi-stable vector
bundles of rank $r$ and trivial determinant on $C$. It is a normal,
projective, unirational variety, of dimension $(r^2-1)(g-1)$. The points of
$\mr$ correspond to isomorphism classes of  vector bundles with trivial
determinant which are  direct sums of  stable vector
bundles of degree zero. 
 The singular locus consists precisely of those  bundles which are
decomposable (with the exception of ${\cal M}_2$ in genus 2, which
 is smooth). The corresponding singularities
are rational Gorenstein -- that is, reasonably mild.

\subsection The Picard group of $\mr$ has been thoroughly studied in [D-N]; let
us recall the main results. Fix some 
$L\in\jg$, and consider the reduced subvariety
$$\Delta _L:=\{E\in\mr\ |\ H^0(E\otimes L)\not=0\}\ .$$
Then $\Delta _L$ is a Cartier divisor in $\mr$; 
the line bundle ${\cal L}:={\cal O}_{\mr}(\Delta _L)$, called the determinant
bundle, is independent of the choice of
$L$ and generates $\Pic(\mr)$. The canonical bundle of $\mr$ is ${\cal
L}^{-2r}$. 

\subsection 
To study the rational map $\varphi _{\cal
L}:\mr\dasharrow|{\cal L}|^*$ associated to the determinant line bundle, the
following construction is crucial. For a vector bundle $E\in\mr$, consider
the locus
$$\Theta _E:=\{L\in \jg\ |\
H^0(C,E\otimes L)\not= 0\}\ .$$Since $\chi (E\otimes L)=0$ for $L$ in
$\jg$, it is readily seen that 
$\Theta _E$ is in a natural way a divisor in 
$\jg$ -- unless it is equal to $\jg$. The latter case (which may occur only for
special bundles) is a serious source of trouble -- see \S 2 below. In the
former case we say that
$E$ admits a theta divisor; this divisor belongs  to the linear system $|r\Theta|
$, where $\Theta
$ is the canonical Theta divisor in $\jg$. 
In this way we get a rational map
$$ \theta: \mr\dasharrow |r\Theta|\ .$$  
\thr Proposition \bib [BNR]
\enonce There is a canonical isomorphism $|{\cal L}|^*\iso |r\Theta |$ which
identifies $\varphi _{\cal L}$ to $\theta $.
\endth
\ind As a consequence, the base locus of the linear system $|{\cal L}|$ is the 
indeterminacy locus of $\theta $, that is, 
 the set of $[E]\in\mr$ such that
$H^0(X,E\otimes L)\not=0$ for all $L\in J^{g-1}$.
\subsection The 
$r$\tx torsion subgroup $J[r]$ of $J$ acts on $\mr$ by tensor product; it also
acts on $|r\Theta |$ by translation, and the map $\theta$ is equivariant with
respect to these actions. In particular, the image of $\theta$ is $J[r]$\tx
invariant.
\ind (1.6) The case when $\theta$ is a morphism is much easier to analyze:
we know then that it is  finite  (since $|{\cal L}|$ is ample, $\theta$ cannot
contract any curve), we know its degree by the Verlinde formula, etc.
Unfortunately there are few cases where this is known to happen:
\th Proposition 
\enonce The base locus of $|{\cal L}|$ is empty in the following cases:
\ind {\rm a)} $r=2\ ;$
\ind {\rm b)} $r=3$, $g=2$ or $3\ ;$
\ind {\rm c)} $r=3$, $C$ is generic.
\endth
\ind All these results except the case $r=g=3$ are due to Raynaud [R]. While a)
 and the first part of b) are easy, c) and the second part of b) are much more
involved. We will discuss the latter in \S 5 below. The proof of c) is reduced,
through a degeneration argument,  to an analogous
statement for torsion-free sheaves on a rational curve with $g$ nodes. 

\section{Base locus}
\subsection Recall that the {\it slope} of a vector bundle $E$ of rank $r$
and degree $d$ is the rational number
$\mu=d/r$. It is convenient to extend the definition of the theta divisor
 to vector bundles $E$
 with {\it integral} slope $\mu$, by putting $\Theta _E:=\{L\in J^{g-1-\mu}\
|\ H^0(C,E\otimes L)\not= 0\}$. If $\delta $ is a line bundle such that $\delta
^{{\scriptscriptstyle\otimes }r}\cong \det E$, the vector bundle
$E_0:=E\otimes \delta ^{-1} $ has trivial determinant and
$\Theta _{E_0}\i \jg$ is the translate by $\delta $ of $\Theta _E\i
J^{g-1-\mu}$.

\subsection  We have the following relations between
stability and existence of the theta divisor:
\indp \numo a)  If $E$ admits a theta divisor, it is semi-stable; 
\indp \numo b) If moreover $\Theta _E$ is a prime divisor, $E$ is stable.
\ind Indeed let $F$ be a proper subbundle of $E$. If 
$\mu(F)>\mu(E)$, the Riemann-Roch theorem implies $H^0(C,F\otimes L)\not=
0$, and therefore $H^0(C,E\otimes L)\not=0$, for all $L$ in $J^{g-1-\mu}$.
 If $\mu(F)=\mu(E)$, one
has $\Theta _E=\Theta _F+\Theta _{E/F}$, so that $\Theta _E$ is not prime.
\subsection The converse of these assertions do not hold. We will see in
(2.6) examples of stable bundles with a reducible theta divisor. The first
examples  of stable bundles with no theta divisor are due to Raynaud [R].
They are restrictions of
projectively flat vector bundles on $J$.
Choose a theta divisor $\Theta $ on $J$. The line
bundle
${\cal O}_J(n\Theta )$ is invariant under the $n$-torsion subgroup $J[n]$ of
$J$. The action of $J[n]$ does {\it not} lift to  ${\cal O}_J(n\Theta )$,
but it does lift to the vector bundle $H^0(J,{\cal O}_J(n\Theta ))^*\otimes
_{\C}{\cal O}_J(n\Theta )$. Thus this vector bundle is the pull back under the
multiplication 
$ n^{}_J:J\rightarrow J$ of a vector bundle $E_n$ on $J$. Restricting $E_n$ to 
the curve 
$C$ embedded in $J$ by an Abel-Jacobi mapping gives the Raynaud bundle
$R_n$. It is well defined up to a twist by an element of $J$, has rank $n^g$ and
slope ${g\over n}$. It has the property that $H^0(C,R_n\otimes \alpha )\not=0$
for all $\alpha \in J$.  Thus if $n\mid g$ $R_n$ has integral slope
and no theta divisor. More generally, Schneider has shown that  a general
vector bundle on
$C$ of rank $n^g$, slope $g-1$ and containing $R_n$ is still stable ([S2]; see
6.6 below). This gives a very large dimension for the base locus of $|{\cal
L}|$,  approximately $(1-{1\over n})\dim\mr$ if
$r=n^g$. Some related results have also been obtained in [A]. 
\subsection Another series of examples have been
constructed by Popa [P].  Let $L$ be a line bundle on $C$ spanned by its global
sections. The {\it evaluation bundle} $E_L$ is defined by the exact sequence 
$$0\rightarrow E_L^*\longrightarrow H^0(L)\otimes _{\C}{\cal O}_C\qfl{ev}
L\rightarrow 0\ ;$$it has the same degree as $L$ and rank $h^0(L)-1$. In
particular, if we choose $\deg L=g+r$ with $r\ge g+2$, $E_L$ has rank $r$ and
slope $\mu=1+{g\over r}$. Then, {\it for all $p$ such that $2\le p\le r-2$ and
$p\mu\in\Z$, the vector bundle $\ext^pE_L$ does not admit a theta divisor}
(see  [S1]). Taking for instance $r=2g$ gives a base point of $|{\cal L}|$ in
${\cal M}_{g(2g-1)}$.
\subsection An interesting limit case of this construction is when $\mu=2$;
this occurs when $L=K_C$, or $r=g$.
The first case has been studied in [FMP]. It turns out that 
the vector bundle $\ext^pE_K$ has a theta divisor, equal to $C_{g-p-1}-C_p$
(here $C_k$ denotes the locus of effective divisor classes in $J^k$). While the
proof is
 elementary for $p=1$, it is extremely involved for the higher exterior powers: it
requires 
going to the moduli space of curves and computing various divisor classes in
the Picard group of this moduli space. 
It remains a challenge to find a direct proof.
\subsection The case $\deg L=2g$ is treated in [B4], building on the
results of [FMP]. Here again $\ext^pE_L$ admits a theta divisor, at least if
$L$ is general enough; it has two components, namely  $C_{g-p-1}-C_p$ and the
translate of  $C_{g-p}-C_{p-1}$ by the class $[K\otimes L^{-1} ]$. These are the
first examples defined on a general curve of stable bundles with a reducible
theta divisor.
\ind (2.7)  Since $|{\cal L}|$ has usually a large base locus, it is natural to
look at the systems $|{\cal L}^k|$ to improve the situation. There
has been much progress on this question in the recent years:
\th Proposition 
\enonce {\rm (i) [P-R]} $|{\cal L}^k|$ is base point free on $\mr$
for $k\ge [{r^2\over 4}]$.
\ind {\rm (ii) [E-P]} For $k\ge r^2+r$, the linear system $|{\cal L}^k|$ defines an
injective morphism of $\mr$ into $|{\cal L}^k|^*$, which is an embedding on
the stable locus.
\endth
\ind On the other hand Popa [P] has observed that one should not be too
optimistic, at least if one believes in the {\it strange duality} conjecture (see
[B2]): this conjecture implies that for $n\mid g$ the Raynaud bundle $R_n$,
twisted by an appropriate line bundle, is a base point of $|{\cal L}^k|$ when
$k\le n(1-{n\over g})$.
\section{Rank 2}
\subsection In rank 2 the situation is now well understood. As pointed out in
(1.6),
$\theta: {\cal M}_2\rightarrow |2\Theta |$ is a finite morphism. In genus 2, 
$\theta $ is actually an {\it isomorphism} onto $\P^3$ [N-R1]. If $X$ is
hyperelliptic of  genus $g\ge 3$, it follows from [D-R] and [B1] that $\theta $
factors through the involution $\iota ^*$ induced by the hyperelliptic
involution and embeds
${\cal M}_{2}/\langle\iota ^*\rangle$ into $|2\Theta |$; moreover the image
admits an explicit geometric description [D-R]. 
\ind (3.2) In the non-hyperelliptic case, after much effort we have now a
complete answer, which is certainly one of the highlights of the subject:
\th{Theorem} 
\enonce If $X$ is not hyperelliptic, $\theta
:{\cal M}_{2} \mono |2\Theta |$ is an embedding.
\endth
\ind The fact that $\theta $ embeds the stable locus of ${\cal M}_{2}$ is 
proved in [B-V], and the remaining part in [vG-I]. Both parts are highly 
nontrivial, and involve some beautiful geometric constructions.

\subsection Thus we can identify ${\cal M}_2$ with a subvariety of $|2\Theta
|\cong
\P^{2^g-1}$, canonically associated to $C$, of   dimension $3g-3$ (1.1).
This variety is invariant under the natural action of $J[2]$ on $|2\Theta |$
(1.5).  Its degree can be computed from the Verlinde formula (see e.g. [Z],
Thm. 1$(iii)$): 
$$\deg {\cal M}_2=(3g-3)!\,2^g(2\pi )^{2-2g}\zeta (2g-2)\ ,$$which gives 
$\deg {\cal M}_2=1$ for $g=2$,  $4$ for $g=3$, $96$ for $g=4$, etc. 

\ind  The singular locus 
$\Sing{\cal M}_{2}$ is
the locus of decomposable bundles in ${\cal M}_{2}$ (1.1), which are 
of the form $\alpha\oplus\alpha ^{-1} $, for 
$\alpha \in J$; the map $\alpha \mapsto \alpha \oplus \alpha
^{-1}$ identifies 
$\Sing{\cal M}_{2}$ to the {\it Kummer variety} ${\cal K}$of $J$ -- that is, the
quotient of
$J$ by the involution $\alpha \mapsto \alpha ^{-1} $. The restriction of
$\theta$ to ${\cal K}=\Sing{\cal M}_2$ is the classical embedding of ${\cal K}$
in $|2\Theta |$, deduced from the map $\alpha \mapsto \Theta _\alpha
+\Theta_{-\alpha}$ from $J$ to $|2\Theta |$. 

\subsection  The case $g=3$, which had been treated previously in
[N-R2], is particularly interesting: we obtain a hypersurface in $|2\Theta |$, of
degree 4, which is $J[2]$\tx invariant and singular along the Kummer variety.
Now Coble shows in [C2] that there is a unique such quartic  (the $J[2]$\tx
invariance is actually superfluous, see [B5]).
Thus in genus 3, the theta map identifies ${\cal M}_{2}$ with the Coble quartic
hypersurface.
\ind Coble gives an explicit equation for this hypersurface, which we now
express in modern terms. Recall that Mumford's theory of the Heisenberg group
allows to find canonical coordinates
$(X_v)_{v\in V}$in the projective space $|2\Theta |$;
here $V$ is a $3$-dimensional subspace of the vector space $J[2]$
over $\F_2$. Then Coble equation reads:
$$\alpha \sum_{u\in V}X_u^4\ + \sum_{\ell=\{u,v\}}\alpha_{d(\ell
)}\,X_u^2\,X_v^2+\sum_{P=\{t,u,v,w\}}\alpha_{d(P)}X_tX_uX_vX_w$$
where the second sum (resp. the third) is taken over the set of affine
lines (resp. planes) in $V$, and $d(\ell )\in\P(V)$ (resp. $d(P)\in
\P(V^*)$) denotes the direction of the line $\ell $ (resp. of the plane
$P$).

\ind In many ways the Coble quartic  ${\cal Q}\subset\P^7$ can be seen as an
analogue of the Kummer quartic surface in $\P^3$. Pauly has proved that ${\cal
Q}$ shares a  famous property of the Kummer surface, the {\it self-duality} : the
dual hypersurface  ${\cal Q}^*\subset(\P^7)^*$ is isomorphic to ${\cal Q}$ [Pa].
The proof is geometric, and includes several beautiful geometric constructions
along the way.
\subsection In genus 4, ${\cal M}_{2}$ is a variety of dimension 9 and degree 96
in
$\P^{15}$.  Oxbury and Pauly have observed that there exists a unique
$J[2]$\tx invariant quartic hypersurface singular along ${\cal M}_{2}$ [O-P]. A
geometric interpretation of this quartic is not known.

\ind (3.6) In arbitrary genus, the quartic hypersurfaces in
$|2\Theta |$ containing ${\cal M}_2$ have been studied in [vG] and [vG-P].
Here is one sample of their results:
\th Proposition  
\enonce Assume that $C$ has no vanishing thetanull. A $J[2]$\tx
invariant quartic form  $F$ on $|2\Theta |$ vanishes on ${\cal
M}_2$ if and only if the hypersurface $F=0$ is singular along ${\cal K}$.
\endth
\ind (Note that though the action of $J[2]$ on $|2\Theta |$ does not come from
a linear action, it {\it does} lift to the space of
quartic forms on $|2\Theta |$. Requiring the invariance of $F$ is stronger than
the invariance of the corresponding hypersurface.)
\ind Van Geemen and Previato also describe
the quartics containing ${\cal M}_2$ in terms of  the Prym varieties associated
to $C$ -- this is related to the Schottky-Jung configuration studied by Mumford. 
  
\section{Genus 2}
\subsection Going to higher rank, it is natural to look first at the genus 2 case.
There a curious numerical coincidence occurs, namely
$$\dim \mr=\dim |r\Theta| = r^2-1\ . $$ 
\ind Recall  that $\theta $ is a finite morphism for
$r=2,3$ (1.6). However already for $r=4$ it is only a rational map: the
Raynaud bundle $R_2$ has rank 4 and slope 1 (2.3), so once twisted
by  appropriate line bundles of degree $-1$ it provides finitely many (actually
16) base points of $|{\cal L}|$.
\ind We have seen that $\theta $ is an isomorphism in rank 2.  In rank 3 there
is again a beautiful story, surprisingly analogous to the rank 2, genus 3 case.
Indeed the Coble quartic has a companion, the {\it Coble cubic} :
this is the unique  cubic hypersurface ${\cal C}\subset
|3\Theta |^*$ singular along $J^2$ embedded in  $ |3\Theta |^*$  by 
the linear system $|3\Theta |$ (this is implicit in Coble [C1]; see [B5] for a
modern explanation).   
\th Theorem
\enonce The map $\theta :{\cal M}_{3}\rightarrow |3\Theta |$ is a double 
covering; the corresponding involution of ${\cal M}_{3}$ is $E\mapsto \iota
^*E^*$, where
$\iota $ is the hyperelliptic involution. The branch locus 
${\cal S}\subset |3\Theta |$ of $\theta $ is a sextic hypersurface, which is the
dual of the Coble cubic ${\cal C}\subset |3\Theta |^*$.
\endth
\ind This is fairly straightforward (see [O]) except for the duality statement,
which  was conjectured by Dolgachev  and proved  in [O] (a different
proof appears in [N]).

\subsection Like for the Coble
quartic  we get an explicit equation for ${\cal C}$ by choosing   a level 3 
structure on $C$, which provides canonical coordinates $(X_v)_{v\in V}$ on
$|3\Theta |^*$, where $V$ is a 2-dimensional
vector space over $\F_3$. Then from [C1] we get the following equation for
${\cal C}$:
$$\alpha_0\sum_{ v\in V} X_v^3\
+6\sum_{\ell =\{u,v,w\}} \alpha_{d(\ell )} X_u
X_v X_w=0\ ,$$where the second sum is taken over the set of affine
lines in $V$, and $d(\ell )\in \P(V)$ is the direction of the
line $\ell $. The 5 coefficients $(\alpha _i)$ satisfy the Burkhardt equation
$$\alpha_0^4-\alpha_0 \sum_{p\in\P(V)}\alpha
_p^3+3\prod_{p\in\P(V)}\alpha _p=0
$$(see [H], 5.3).

\subsection In rank $r\ge 4$ we start getting base points, and this causes a lot of
trouble -- since $\theta$ is only rational, we cannot compute its degree using
intersection theory. However we still have:
\thr Proposition \bib [B6]
\enonce The rational map $\theta :\mr\dasharrow |r\Theta |$ is generically
finite {\rm (or, equivalently, dominant).}
\endth
\ind The idea is to prove the finiteness of $\theta ^{-1} (\Theta +\Delta )$, 
where $\Delta
$ is a general element of $|(r-1)\Theta| $. Any decomposable bundle in that fibre
must be of the form ${\cal O}_C\oplus F$ for some $F\in{\cal M}_{r-1}$ with
$\Theta _F=\Delta $; reasoning by  induction on $r$ we can assume that there
are finitely many such $F$.  Thus the whole point is to control the {\it stable}
bundles
$E$ with
$\Theta _E=\Theta +\Delta $. Now the condition $\Theta _E\supset \Theta $
means by definition
$H^0(C,E(p))\not=0$ for all $p\in C$, or equivalently
$H^0(C,E'(-p))\not=0$ for all $p\in C$, where $E':=E^*\otimes K_C^{-1} $ is the
Serre dual of $E$. Since $h^0(E')=r$ by stability of $E$, this implies that
the global sections of $E'$ generate a subbundle of rank $<r$. A precise analysis
of this situation allows to prove
that there are only finitely many such bundles $E$ with $\Theta _E=\Theta
+\Delta $.

\subsection The map $\theta$ is no longer finite in rank $r\ge 4$, in fact 
 it admits some fibres of dimension $\ge [{r\over 2}]-1$ [B6]. When
$r$ is even, this is seen by restricting $\theta$ to the moduli space of {\it
symplectic} bundles: the corresponding moduli space has dimension ${1\over
2}r(r+1)$, but its image under $\theta$ lands in the subspace $|r\Theta |^+$ of
$|r\Theta |$ corresponding to even theta functions of order $r$, which has
dimension ${r^2\over 2}+1$. For $r$ odd one considers bundles of the form
${\cal O}_C\oplus F$ with $F$ symplectic.

 \subsection It would be interesting to find the {\it degree} of $\theta$, which is
unknown already in genus 4. Note that for trivial reasons it has to be 
exponential in $g$ (see [B6], 2.3). Brivio and Verra have found a nice geometric
interpretation of the generic fibre of $\theta$  which  might lead at least to a
good estimate for 
$\deg \theta $ [V].

\section{Higher rank and genus}  
\ind Not much is known here. We already mentioned the following result
proved in [B6]:
\th Proposition
\enonce In genus $3$ the map $\theta:{\cal M}_{3}\rightarrow |3\Theta |$ is a 
finite morphism.
\endth
\ind The proof is rather roundabout, and gives actually a more interesting result:
the complete list of stable vector bundles $E$ of rank 3 and degree 0 such that
$\Theta _E\supset \Theta $. It turns out that the bundles in this list have a theta
divisor, so Proposition 5.1 follows. 
\subsection The idea for establishing that list is to translate the problem into a
classical question of projective geometry. Similarly to the genus 2 case, the
condition
$\Theta _E\supset \Theta $ means $H^0(E(p+q))\not= 0$ for all $p,q$ in $C$, 
or equivalently $H^0(E'(-p-q))\not= 0$, where $E':=E^*\otimes K_C^{-1} $ is the
Serre dual of $E$. One checks that stability implies $h^0(E')=6$ and
$h^0(E'(-p))=3$ for $p$ general in $C$. This gives a family of 2-planes in
$\P(H^0(E))\cong\P^5$, parametrized by $C$, such that any two planes of the
family intersect. It turns out that the maximal such families have been classified
in a beatiful paper by Morin [M]: there are three families given by linear algebra
(like the 2-planes contained in a given hyperplane), and three coming from
geometry: the 2-planes contained in a smooth quadric, the tangent planes to the
Veronese surface, and the planes intersecting the Veronese surface along a conic.
Translating back this result in terms of vector bundles gives the list we
were looking for.

\subsection This list also shows that $\theta^{-1} (\Theta +\Theta_F)=\{{\cal
O}_C\oplus F\}$ for $F$ general in ${\cal M}_{2}$. This might indicate that
$\theta$ has degree one; this would follow if we could prove the injectivity
of its tangent map at ${\cal O}_C\oplus E$ for some
$E$ in ${\cal M}_2$, perhaps in the spirit of [vG-I].

\section{Questions and conjectures}
\ind The list of results ends at this point, but let me finish with a (small) list of
open problems. About the general behaviour of the theta map, the most
optimistic statement would be:
\th Speculation 
\enonce For $g\ge 3$, $\theta$ is generically injective if
$C$ is not hyperelliptic, and generically two-to-one onto its image if $C$ is
hyperelliptic.
\endth  
\ind Note that 
in the hyperelliptic case $\theta$ factors as in thm.\ 4.2 through the
non-trivial involution $E\mapsto \iota^*E^*$. 
 Admittedly the evidence for 6.1 is very weak: the only
case where it is known is in rank 2.\smallskip 
\ind As for base points, Proposition 1.6 leads naturally to:
\th Conjecture
\enonce Every  bundle $E\in{\cal M}_3$ has a theta divisor.
\endth
\subsection There is an integer $r(C)$ such that $\theta $ is a
morphism for $r<r(C)$ but only a rational map for $r\ge r(C)$ (observe that if
$E\in\mr$ has no theta divisor, so does $E\oplus F$ for any $F$ in ${\cal
M}_s,\ s\ge 1$). We know very little about this integer: we have
$r(C)=4$ for
$g=2$, $4\le r(C)\le 8$ for $g=3$,  and $r(C)\le {1\over 2}(g+1)(g+2)$ [A].
\th Questions
\enonce {\rm a)} Does $r(C)$ depend only on $g${\rm ?} 
\ind {\rm b)} Put
$r(g):=\min r(C)$ for all curves $C$ of genus $g$. Is $r(g)$
 an increasing function of $g${\rm ?}
\endth
\ind The next question do not involve directly the theta map, but 
it is related to several  questions about the existence of  theta divisors. 

\th Conjecture 
\enonce Let $\pi :X\rightarrow Y$ be a finite morphism between smooth
projective curves of genus $\ge 2$. The direct image
$\pi _*L$ of a general vector bundle 
$L$ on $X$ is stable.
\endth
\ind One reduces readily to the case when $L$ is a line bundle. The problem
depends on a crucial way on  the degree of $L$:  one can prove for instance 
that $\pi _*L$ is stable (for $L$ generic) if $|\chi (L)|<g+{g^2\over r}$, where
$r$ is the degree of $\pi $ and $g$ the genus of $Y$ (see [B3]).
\ind One of the relations between this conjecture and
the existence of  theta divisors is the following: the conjecture  for $L$ 
general of degree $d$ is implied by the existence of {\it a vector bundle $E$ of
rank $r$ and degree
$g(X)-1-d$ such that
$\pi ^*E$ admits a prime theta divisor}. Indeed  
we have $\Theta _{\pi _*L\otimes E}=(\pi ^*)^{-1}( \Theta _{L\otimes \pi
^*E})$; if $\Theta _{\pi^*E}$ is prime, so is $\Theta _{\pi _*L\otimes E}$ for
general $L$, and as in (2.2) this implies that $\pi _*L$ is stable.

\vskip2cm 
\centerline{ REFERENCES} \vglue15pt\baselineskip12.8pt
\def\num#1{\smallskip\item{\hbox to\parindent{\enskip [#1]\hfill}}}
\parindent=1.3cm 
\def\pp#1{Preprint {\tt math.AG/{#1}}.}
\num{A} D. {\pc ARCARA}: {\sl A lower bound for the dimension of the base 
locus of the generalized theta divisor}. C. R. Math. Acad. Sci. Paris {\bf 340}
(2005), 131--134.
\num{B1} A. {\pc BEAUVILLE}:  {\sl 	Fibr\'es de rang $2$ sur les courbes,
fibr\'e d\'eterminant et fonctions th\^eta}. Bull. Soc. math. France {\bf 116} 
(1988), 431--448. 
\num{B2} A. {\pc BEAUVILLE}: {\sl Vector bundles on curves and generalized
theta functions: recent results and open problems}. Current Topics in Complex
Algebraic Geometry, MSRI Publications {\bf 28}, 17-33; Cambridge University
Press (1995).
\num{B3}  A. {\pc BEAUVILLE}: {\sl On the stability of the direct image of a
generic vector bundle}. Preprint (2000), available at {\tt
http://math1.unice.fr/$\!\sim$beauvill/ pubs/imdir.pdf}
\num{B4}  A. {\pc BEAUVILLE}: {\sl Some stable vector bundles  with reducible
theta divisors}. Manuscripta Math. {\bf 110} (2003), 343--349.
\num{B5} A. {\pc BEAUVILLE}: {\sl The Coble hypersurfaces}. C. R. Math. 
Acad. Sci. Paris {\bf 337} (2003),  189--194.

\num{B6} A. {\pc BEAUVILLE}: {\sl Vector bundles and theta
functions on curves of genus $2$ and $3$}. \pp{0406030}  To
appear.
\num{BNR} A. {\pc BEAUVILLE}, M.S. {\pc NARASIMHAN}, S. {\pc
RAMANAN}: {\sl Spectral curves and the generalised theta
divisor}. J. Reine Angew. Math. {\bf 398} (1989), 169--179. 

\num{B-V} S. {\pc BRIVIO}, A. {\pc VERRA}: {\sl The theta divisor of  ${\rm
SU}\sb C(2,2d)\sp s$ is very ample if $C$ is not hyperelliptic}. Duke Math. J.
{\bf 82} (1996),  503--552.
 \num{C1} A. {\pc COBLE}: {\sl Point Sets and Allied Cremona
Groups} III. Trans.  Amer. Math. Soc. {\bf  18} (1917),  331--372. 
\num{C2} A. {\pc COBLE}: {\sl Algebraic geometry and theta functions}.
Amer. Math. Soc. Colloquium Publi. 10 (1929). Amer. Math. Soc.,
Providence (1982).

\num{D-R} U.V. {\pc DESALE}, S. {\pc RAMANAN}: {\sl Classification of vector
bundles of rank $2$ on hyperelliptic curves}. Invent. math. {\bf 38} (1976),
161--185.

\num{D-N} J.M. {\pc DREZET}, M.S. {\pc NARASIMHAN}: {\sl Groupe de Picard
des vari\'et\'es de modules de fibr\'es semi-stables sur les courbes
alg\'ebriques.} Invent. math. {\bf 97} (1989), 53--94.

\num{E-P} E. {\pc ESTEVES}, M. {\pc POPA}: {\sl Effective very ampleness for
generalized theta divisors}. Duke Math. J. {\bf 123} (2004),  429--444.
\num{vG} {\pc B. VAN} {\pc GEEMEN}:
{\sl Schottky-Jung relations and 
vector bundles on hyperelliptic curves}. Math. Ann.  {\bf 281} 
(1988),  431--449.
\num{vG-I} {\pc B. VAN} {\pc GEEMEN}, E. {\pc IZADI}: {\sl The tangent
space to the moduli space of vector bundles on a curve and the singular
locus of the theta divisor of the Jacobian}. J. Algebraic Geom. {\bf 10}
(2001),  133--177. 

\num{vG-P} {\pc B. VAN} {\pc GEEMEN}, E. {\pc PREVIATO}: {\sl Prym varieties
and the Verlinde formula.} Math. Annalen {\bf 294} (1992), 741--754.

\num{FMP} G. {\pc FARKAS}, M. {\pc MUSTA\c{T}\v{A}}, M. {\pc POPA}:
{\sl Divisors  on ${\cal M}_{g,g+1}$ and the Mini\-mal Resolution Conjecture
for points on canonical curves}. Ann. Sci. \'Ecole Norm. Sup. (4) {\bf 36} (2003),
\no 4, 553--581.

\num{H} B. {\pc HUNT}: {\sl The geometry of some special arithmetic 
quotients}. Lecture Notes in Math. {\bf 1637}. Springer-Verlag, Berlin
(1996).
\num{K-N} S. {\pc KUMAR}, M.S. {\pc NARASIMHAN}:  
{\sl Picard group of the moduli spaces of $G$-bundles}. 
Math. Ann. {\bf 308} (1997),  155--173.

\num{M} U. {\pc MORIN}: {\sl Sui sistemi di piani a due a due incidenti}. Atti
Ist. Veneto {\bf 89} (1930), 907--926.

\num{N-R1} M.S. {\pc NARASIMHAN}, S. {\pc
RAMANAN}: {\sl Moduli of vector bundles on a compact Riemann
surface}. Ann. of Math. (2) {\bf 89} (1969), 14--51. 

\num{N-R2} M.S. {\pc NARASIMHAN}, S. {\pc RAMANAN}: {\sl
$2\theta$\tx linear  systems  on Abelian varieties.} Vector bundles
on algebraic varieties, 415--427; Oxford University Press (1987).
\num{N} Q.M. {\pc NGUYEN}: {\sl The moduli space of rank $3$ vector bundles 
with trivial determinant over a curve of genus $2$ and duality}. \pp{0408318}

\num{O} A. {\pc ORTEGA}:  {\sl On the moduli space of rank 3 vector
bundles on a genus 2 curve and the Coble cubic}. J. Algebraic Geom. {\bf 
14} (2005), 327--356.
\num{O-P} W. {\pc OXBURY}, C. {\pc PAULY}: {\sl Heisenberg invariant 
quartics and ${\cal SU}_C(2)$ for a curve of
   genus four}. Math. Proc. Cambridge Philos. Soc. {\bf  125} (1999),
295--319. 

\num{Pa} C. {\pc PAULY}: {\sl Self-duality of Coble's quartic hypersurface and
applications}. Michigan Math. J. {\bf 50} (2002),  551--574.

\num{P} M. {\pc POPA}: {\sl On the base locus of the generalized theta divisor}.
C. R. Acad. Sci. Paris S\'er. I Math. {\bf 329} (1999),  507--512.
\num{P-R} M. {\pc POPA}, M. {\pc ROTH}: {\sl Stable maps and Quot schemes}.
Invent. Math. {\bf 152} (2003),  625--663.

\num{R} M. {\pc RAYNAUD}: {\sl Sections des fibr\'es vectoriels sur
une courbe}. Bull. Soc. Math. France {\bf 110} (1982),
103--125. 
\num{S1} O. {\pc SCHNEIDER}: {\sl Stabilit\'e des fibr\'es $\ext^{p}E_{L}$
et condition de Raynaud}. Preprint {\tt math.AG/0309277}, to appear in Ann.
Fac. Sci. Toulouse Math.
\num{S2} O. {\pc SCHNEIDER}: {\sl Sur la dimension de l'ensemble des points 
base du fibr\'e d\'eterminant sur
  l'espace des modules des fibr\'es vectoriels semi-stables sur une courbe}.
Preprint {\tt math.AG/0501320}.
\num{V} A. {\pc VERRA}: Lecture at the Conference ``Algebraic cycles and
motives", Leiden (2004).
\num{Z} D. {\pc ZAGIER}: {\sl Elementary aspects of the Verlinde formula and of
the Harder-Narasimhan-Atiyah-Bott formula}. Proceedings of the Hirzebruch 65
Conference on Algebraic Geometry, 445--462, Israel Math.
Conf. Proc. {\bf 9} (1996).
\vskip1cm
\def\pc#1{\eightrm#1\sixrm}
\hfill\vtop{\eightrm\hbox to 5cm{\hfill Arnaud {\pc BEAUVILLE}\hfill}
 \hbox to 5cm{\hfill Institut Universitaire de France\hfill}\vskip-2pt
\hbox to 5cm{\hfill \&\hfill}\vskip-2pt
 \hbox to 5cm{\hfill Laboratoire J.-A. Dieudonn\'e\hfill}
 \hbox to 5cm{\sixrm\hfill UMR 6621 du CNRS\hfill}
\hbox to 5cm{\hfill {\pc UNIVERSIT\'E DE}  {\pc NICE}\hfill}
\hbox to 5cm{\hfill  Parc Valrose\hfill}
\hbox to 5cm{\hfill F-06108 {\pc NICE} Cedex 2\hfill}}
\end